\documentclass[12pt]{article}
\usepackage{pdfrender,xcolor}
\usepackage{graphicx,color}
\usepackage{amsmath,amsfonts,amssymb,amsthm}
\usepackage[hidelinks,colorlinks=true,linkcolor=blue,citecolor=blue]{hyperref}
\usepackage{rsfso}
\usepackage{comment}
\usepackage{enumerate}
\usepackage{pdfpages}

\usepackage[numbers]{natbib}
\newcommand{\vc}[1]{{\boldsymbol #1}}

\newcommand{\sr}[1]{{\mathcal #1}}
\newcommand{\dd}[1]{\mathbb{#1}}

\newcommand{\qvs}[1]{\big[ #1 \big]}

\newcommand{\sgn}{{\rm sgn}}

\newcommand{\eq}[1]{(\ref{eq:#1})}
\newcommand{\lem}[1]{Lemma~\ref{lem:#1}}

\newcommand{\thr}[1]{Theorem~\ref{thr:#1}}

\newcommand{\cond}[1]{Condition~\ref{cond:#1}}
\newcommand{\rem}[1]{Remark~\ref{rem:#1}}
\newcommand{\fig}[1]{Figure~\ref{fig:#1}}

\newcommand{\sectn}[1]{Section~\ref{sec:#1}}

\newcommand{\thrt}[1]{\ref{thr:#1}}

\newcommand{\pend}{\hfill \thicklines \framebox(6.6,6.6)[l]{}}

\newenvironment{proof*}[1]{\noindent {\sc  #1} \rm}{\pend}
\newenvironment{keyword}[1]{ {\bf #1} \rm}{}
\newenvironment{MSC}[1]{ {\bf #1} \rm}{}

\newtheorem{theorem}{Theorem}[section]
\newtheorem{lemma}{Lemma}[section]

\newtheorem{condition}{Condition}[section]
\newtheorem{remark}{Remark}[section]

\newcommand{\setnewcounter} {
\setcounter{subsection}{0}
\setcounter{equation}{0}
\setcounter{conjecture}{0}
\setcounter{assumption}{0}
\setcounter{question}{0}
\setcounter{definition}{0}
\setcounter{theorem}{0}
\setcounter{corollary}{0}
\setcounter{lemma}{0}
\setcounter{proposition}{0}
\setcounter{remark}{0}
}

\begin{document}
 \title{\Large \bf The stationary distributions of state-dependent diffusions reflected at one and two sides}

%%%%%%%%%%%%%%
%   AUTHORS  %
%%%%%%%%%%%%%%
\author{Masakiyo Miyazawa\footnotemark}
\date{\today}
%\date{}

\maketitle

\begin{abstract}
Consider a one-dimensional diffusion process which has state-dependent drift and deviation and is reflected at the origin, which is called a one-side reflected diffusion or simply reflected diffusion. We are particularly interested in the case that its drift and deviation are discontinuous. We define this reflected diffusion as the solution of a stochastic integral equation, and find conditions for its positive recurrence, We then derive its stationary distribution under these conditions. As a related problem, we also consider the case that it is reflected at two sides, which is called a two-sides reflected diffusion. Its existence, positive recurrence and stationary distribution are similarly studied.

In the literature, these problems are studied through a state-dependent diffusion on the whole line particularly when the drift and deviation are discontinuous. However, the reflected process itself is not defined in such a study. Thus, the stationary distribution has not been fully studied for a general state-dependent reflected diffusion. We aim to fills this insufficiency and to make the stationary distributions of reflected diffusions widely available in application.
\end{abstract}

\keyword{Keywords:}{state-dependent diffusion, discontinuous diffusion coefficient, reflected at two sides, stationary distribution, stochastic integral equation, generalized Ito formula.}

\MSC{MSC Classification:}{60H20, 60J25, 60J55, 60H30, 60K25}

\footnotetext[1]{Department of Information Sciences,
Tokyo University of Science, Noda, Chiba, Japan}

\section{Introduction}
\label{sec:introduction}

We consider a nonnegative real-valued diffusion process $Z(\cdot) \equiv \{Z(t); t \ge 0\}$ which has state-dependent diffusion coefficients, that is, drift $b(x)$ and deviation $\sigma(x) > 0$, and is reflected at the origin. This $Z(\cdot)$ is called a state-dependent one-side reflected diffusion or simply a reflected diffusion. We define this process $Z(\cdot)$ as the solution of the following stochastic integral equation, SIE for short. 
\begin{align}
\label{eq:Z-SIE}
  Z(t) = Z(0) & + \int_{0}^{t} \sigma(Z(s)) dW(s) + \int_{0}^{t} b(Z(s)) ds + Y(t) \ge 0, \qquad t \ge 0,
\end{align}
where $b(x)$ and $\sigma(x) > 0$ are measurable functions of $x \ge 0$, $W(\cdot)\equiv \{W(t); t \ge 0\}$ is the standard Brownian motion, and $Y(\cdot)\equiv \{Y(t); t \ge 0\}$ is a right-continuous and non-deceasing process satisfying that $\int_{0}^{t} 1(Z(s) > 0) dY(s) = 0$ for $t \ge 0$. For the righthand side of \eq{Z-SIE} to be well defined, it is assumed that the following integral exists almost surely.
\begin{align}
\label{eq:Z-ass 1}
  \int_{0}^{t} \left(\sigma^{2}(Z(s)) + |b(Z(s))|\right) ds < \infty, \qquad t > 0.
\end{align}
A sufficient condition for this is that $\sigma^{2}(x)$ and $|b(x)|$ are locally bounded in $x \ge 0$, but these boundedness may not be needed for \eq{Z-ass 1} to hold. It is also notable that no further assumption is needed for $\sigma(x)$ and $b(x)$ except for their measurability and \eq{Z-ass 1} in our results.

We denote the stochastic basis for stochastic processes by $(\Omega,\sr{F},\dd{F},\dd{P})$, where $\dd{F}$ is a right-continuous filtration on probability space $(\Omega,\sr{F},\dd{P})$, that is, $\dd{F} = \{\sr{F}_{t}; t \ge 0\}$ is the sequence of nondecreasing sub-$\sigma$-field $\sr{F}_{t}$ of $\sr{F}$ satisfying that $\sr{F}_{t} = \cap_{s > t} \sr{F}_{s}$. Throughout the paper, we assume that all stochastic processes are right continuous with lefthand limits.

We are interested in the stationary distribution of $Z(\cdot)$ when it may have discontinuous drift and deviation, motivated by an energy saving problem for a single station (see \cite{Miya2025}). We first prove that SIE \eq{Z-SIE} has a weak solution $(Z(\cdot),Y(\cdot),W(\cdot))$ under an appropriate condition. In this proof, we will use a one-dimensional diffusion process without reflection, which is defined as the solution $X(\cdot) = \{X(t); t \ge 0\}$ of the following SIE.
\begin{align}
\label{eq:X-SIE}
  X(t) = X(0) & + \int_{0}^{t} \sigma(X(s)) dW(s) + \int_{0}^{t} b(X(s)) ds, \qquad t \ge 0,
\end{align}
where a similar condition to \eq{Z-ass 1} is assumed for $X(\cdot)$ instead of $Z(\cdot)$ for the righthand side of \eq{X-SIE} to be well defined. We call this $X(\cdot)$ a state-dependent diffusion or simply a diffusion. We briefly introduce the known results on the existence of $X(\cdot)$.

We then show that $Z(\cdot)$ is a strong Markov process, and give necessary and sufficient condition for $Z(\cdot)$ to be recurrent (\lem{Z-X-existence}). Based on these results, we derive the stationary distribution of $Z(\cdot)$ from its stationary measure under this recurrent condition (\thr{stationary}). These results are proved using the generalized Ito formula and local time. For them, we refine ideas used in \citet{Miya2024b}, which studies a multi-level semi-martingale reflected Brownian motion (SRBM) on $[0,\infty)$, which is a special case of the reflected diffusion such that that $b(x)$ and $\sigma(x)$ are simple functions.

As a related problem, we are also interested in the stationary distribution of the state-dependent diffusion which is reflected at two sides. This process is called a two-sides reflected diffusion. Note that this reflected process may be used for the diffusion approximation of a single station queue with a finite buffer in heavy traffic, so it is interested in application. We prove similar results (\lem{U-X-existence} and \thr{stationary-a}) to those for the one side reflected diffusion (\lem{Z-X-existence} and \thr{stationary}).

In the literature, the stationary distributions of these state-dependent reflected diffusions have not been well studied except for the recent paper \cite{Miya2024b} when the diffusion coefficients are discontinuous. So far, in the literature, a typical approach for the stationary distribution of the reflected diffusions firstly focuses on the diffusion process $X(\cdot)$. Then, derive the stationary measure  of $X(\cdot)$ under appropriate conditions, and restrict it on the half line or a closed interval. This restricted measure is considered as the stationary measure of the corresponding reflected diffusions (see Chapter 23 of \cite{Kall2001}). However, the reflected processes themselves are not defined, so it is unclear for what reflected processes the restricted measures are stationary although the results except for the two-sides reflected diffusion exactly correspond with our results for $Z(\cdot)$ (see \thr{stationary}). This type of approach is used even for the case that the diffusion coefficients are smooth (e.g., see Chapter V of \cite{BhatWaym2009}). We discuss their details in \sectn{concluding}.

This paper is made up by five sections. In \sectn{preliminary}, basic conditions are introduced, and preliminary lemmas are presented for the diffusion and its reflection. In \sectn{stationary}, a main result, \thr{stationary} is presented and proved. In \sectn{extension}, we extend these results for the two-sides reflected diffusion. Finally, in \sectn{concluding}, we discuss the literature about the one and two sides reflected diffusions and their stationary distributions. Our key tools are the generalized Ito formula and local time. They are briefly introduced in the appendix.

\section{Preliminary lemmas}
\label{sec:preliminary}
\setnewcounter

In this section, we prepare basic results for the reflected diffusion process $Z(\cdot)$, and prove its existence  using a slight modification of the diffusion process $X(\cdot)$. For this, we first discuss about the existence of $X(\cdot)$ as a solution of \eq{X-SIE}. Because of the discontinuity of the diffusion coefficients, we cannot assume that they are Lipschitz continuous, which is standard to consider the solution of a stochastic integral equation. As is well known, \citet{EngeSchm1985} found a minimal condition without any continuity on the diffusion coefficients under which SIE \eq{X-SIE} has a unique weak solution when drift $b(x)$ vanishes for all $x \in \dd{R}$, where $\dd{R}$ is the set of all real numbers. This result had been extended to the non-vanishing drift. To present this extension, we define an explosion time $S$ for $X(\cdot)$ as
\begin{align*}
   S = \lim_{n \to \infty} \inf\{t \ge 0; |X(t)| \ge n\}
\end{align*}
Then, the following result can be found in \citet{KaraShre1998} (see also Theorems 21.11 and 23.1 of \citet{Kall2001}).

\begin{lemma}[Theorem 5.15 of \cite{KaraShre1998}]
\label{lem:X-existence}
SIE \eq{X-SIE} has a unique weak solution $(X(\cdot), W(\cdot))$ up to the explosion time $S$ of $X(\cdot)$ if \cond{X-existence} below holds. Furthermore, this $\{X(t); t \in [0,S)\}$ is a strong Markov process.
\end{lemma}
\begin{condition}
\label{cond:X-existence}
The functions $b(x)$ and $\sigma(x)$ satisfy
\begin{align}
\label{eq:sigma positive}
 & \sigma(x) > 0, \qquad \forall x \in \dd{R},\\
\label{eq:local integrability 1}
 & \int_{x_{1}}^{x_{2}} \frac 1{\sigma^{2}(y)} dy < \infty, \qquad \int_{x_{1}}^{x_{2}} |\beta(y)| dy < \infty, \qquad \forall x_{1}, x_{2} \in \dd{R}, x_{1} < x_{2},
\end{align}
where 
\begin{align}
\label{eq:beta}
  \beta(x) = \frac {2b(x)} {\sigma^{2}(x)}, \qquad x \ge 0.
\end{align}
\end{condition}

We now consider the solution $Z(\cdot)$ of SIE \eq{Z-SIE}. Since this SIE does not need $b(x)$ and $\sigma^{2}(x)$ for $x < 0$, we weaken \cond{X-existence} to
\begin{condition}
\label{cond:Z-existence}
The functions $b(x)$ and $\sigma(x)$ satisfy
\begin{align}
\label{eq:sigma positive 2}
 & \sigma(x) > 0, \qquad \forall x \in \dd{R}_{+} \equiv [0,\infty),\\
\label{eq:local integrability 2}
 & \int_{x_{1}}^{x_{2}} \frac 1{\sigma^{2}(y)} dy < \infty, \qquad \forall (x_{1},x_{2}) \in \dd{R}_{+}^{2} \mbox{ satisfying } x_{1} < x_{2},\\
\label{eq:local integrability 3}
 & \int_{0}^{x+\varepsilon} |\beta(y)| dy < \infty, \qquad \forall x \in \dd{R}_{+}, \exists \varepsilon > 0.
\end{align}
\end{condition}

\lem{Z-X-existence} below is a starting point of our study, in which the weak solution of SIE \eq{Z-SIE} is constructed. A key idea of this construction is to define a diffusion on $\dd{R}$ by symmetrically modifying the drift $b(x)$ and deviation $\sigma(x)$ of $X(\cdot)$ for $x < 0$, which does not influence SIE \eq{Z-SIE} as remarked above. Then, the absolute value of this diffusion will be the weak solution of \eq{Z-SIE}. This idea is originally due to \citet{AtarCastReim2024}.

We denote the conditional expectation given $Z(0)=x$ by $\dd{E}_{x}$, and define stopping time $\tau_{y}$ as
\begin{align*}
  \tau_{y} = \inf\{t > 0; Z(t) = y\}, \qquad y \ge 0,
\end{align*}
where $\tau_{y} = \infty$ is $Z(t) \not= y$ for all $t \ge 0$. Then, we call $Z(\cdot)$ is recurrent if
\begin{align}
\label{eq:recurrence}
  \dd{P}_{x}[\tau_{y} < \infty] = 1, \qquad \forall x,y \ge 0. 
\end{align}
Otherwise, we call it transient. Note that \eq{recurrence} implies that $Z(\cdot)$ is Harris recurrent, by which $Z(\cdot)$ has a stationary measure (e.g., see \cite{MeynTwee1993a,MiyaMoro2023} for the definition of Harris recurrence).
\begin{lemma}
\label{lem:Z-X-existence}
Under \cond{Z-existence}, (i) the stochastic integral equation \eq{Z-SIE} has a weak solution $({Z}(\cdot),{Y}(\cdot),{W}(\cdot))$ up to the explosion time of $Z(\cdot)$, and (ii) this ${Z}(\cdot)$ is a strong Markov process up to its explosion time, where the same notations are used for the weak solution for simplicity. Furthermore, define $\eta(x)$, called a scale function, as
\begin{align}
\label{eq:scale-f}
  \eta(x) = \int_{0}^{x} \exp\left(- \int_{0}^{y} \beta(v) dv\right) dy, \qquad x \ge 0,
\end{align}
then (iii) the ${Z}(\cdot)$ is recurrent if and only if
\begin{align}
\label{eq:eta-infty}
  \eta(\infty) \equiv \lim_{x \to \infty} \eta(x) = \infty.
\end{align}

\end{lemma}

\begin{remark}
\label{rem:Z-X-existence}
(i) is proved in Remark 4.1 of \cite{Miya2024b} (see also Lemma 4.3 of \cite{AtarCastReim2024}). We here prove it in a slightly different way in order to prove (ii). (iii) and its proof are similar to those of Corollary 9.3 in Chapter V of \cite{BhatWaym2009}, which considers the diffusion process $X(\cdot)$ on an open interval, assuming that $b(x)$ and $\sigma(x)$ are continuous in the sense of (1.2)' of \cite{BhatWaym2009}. However, these continuities do not hold in our case. Furthermore, the reflected diffusion $Z(\cdot)$ is not considered there. This is further discussed in \sectn{concluding}.
\end{remark}

\begin{remark}
\label{rem:Z-uniqueness}
The uniqueness of the weak solution ${Z}(\cdot)$ is likely true, but has not yet been proved while there are sufficient conditions for it (see Lemma 4.1 of \cite{AtarCastReim2024a} and Lemma 2.2 of \cite{Miya2024b}). On the other hand, \citet{Tana1979} shows that the SIE \eq{Z-SIE} has a unique strong solution if the diffusion coefficients are sufficiently smooth and locally bounded. See \sectn{concluding} for its details.
\end{remark}

\begin{proof}[Proof of \lem{Z-X-existence}]

We prove the part (i). For the given $(X(\cdot),W(\cdot))$, define
\begin{align*}
%\label{eq:sgn-VD}
 & \widetilde{\sigma}(x) = \sigma(|x|), \qquad \widetilde{b}(x) = \sgn(x) b(|x|), \qquad x \in \dd{R}.
\end{align*}
where $\sgn(x) = 1(x >0) - 1(x < 0)$. Then, consider the following stochastic integral equation for $\widetilde{X}(\cdot) \equiv \{\widetilde{X}(t); t \ge 0\}$.
\begin{align}
\label{eq:tX-SIE}
  \widetilde{X}(t) = \widetilde{X}(0) + \int_{0}^{t} \widetilde{\sigma}(\widetilde{X}(s)) dW(s) + \int_{0}^{t} \widetilde{b}(\widetilde{X}(s)) ds, \qquad t \ge 0.
\end{align}

Since $b(x)$ and $\sigma(x)$ satisfy \cond{Z-existence}, it is easy to see that $\widetilde{\sigma}(x)$ and $\widetilde{b}(x)$ satisfy \cond{X-existence} instead of $\sigma(x)$ and $b(x)$. Hence, by \lem{X-existence}, the SIE \eq{tX-SIE} has the weak solution $\widetilde{X}(\cdot)$ which is unique in the weak sense and a strong Markov process up to its explosion time. On the other hand, from Tanaka formula \eq{Tanaka} for $a = 0$, we have
\begin{align}
\label{eq:|tX|-SIE}
 & |\widetilde{X}|(t) = |\widetilde{X}|(0) + \int_{0}^{t} \sgn (\widetilde{X}(s)) \left(\widetilde{\sigma}(\widetilde{X}(s)) dW(s) + \widetilde{b}(\widetilde{X}(s)) ds\right)+ \widetilde{L}_{0}(t) \nonumber\\
 & \quad = \int_{0}^{t}  \left(\sigma(|\widetilde{X}|(s)) \sgn(\widetilde{X}(s)) dW(s) + b(|\widetilde{X}|(s)) ds\right)+ \widetilde{L}_{0}(t),
\end{align}
because $\sgn^{2}(\widetilde{X}(s)) = 1$, where $|\widetilde{X}|(t) = |\widetilde{X}(t)|$, and  $\widetilde{L}_{0}(\cdot) \equiv \{\widetilde{L}_{0}(t); t \ge 0\}$ is the local time of $\widetilde{X}(\cdot)$ at $0$. Define
\begin{align*}
  \widetilde{Z}(t) = |\widetilde{X}|(t), \qquad \widetilde{Y}(t) = \widetilde{L}_{0}(t), \qquad \widetilde{W}(t) = \int_{0}^{t}\sgn(\widetilde{X}(s)) dW(s),
\end{align*}
then $\widetilde{W}(\cdot) \equiv \{\widetilde{W}(t); t \ge 0\}$ is the standard Brownian motion and
\begin{align}
\label{eq:substitution-f}
  \int_{0}^{t} \sigma(|\widetilde{X}|(s)) \sgn(\widetilde{X}(s)) dW(s) = \int_{0}^{t}  \sigma(|\widetilde{X}|(s)) d\widetilde{W}(s),
\end{align}
by the extended version of Theorem 2.12 of \cite{ChunWill1990} (see its Section 3.3). Hence, \eq{|tX|-SIE} shows that $(\widetilde{Z}(\cdot), \widetilde{Y}(\cdot), \widetilde{W}(\cdot))$ is the weak solution of \eq{Z-SIE} up to the explosion time of $\widetilde{Z}(\cdot)$. Thus, rewriting this weak solution as $(Z(\cdot), Y(\cdot), W(\cdot))$ by our convention, we have proved (i).

For the parts (ii) and (iii), we use the notations $\widetilde{X}(\cdot)$, $\widetilde{W}(\cdot)$, $\widetilde{Z}(\cdot)$ and $\widetilde{Y}(\cdot)$ which are defined in the part (i). We now prove (ii). Let $\widetilde{\dd{F}} \equiv \{\widetilde{\sr{F}}_{t}; t \ge 0\}$ be the natural right-continuous filtration generated by $\widetilde{X}(\cdot)$ and $\widetilde{W}(\cdot)$, let $\tau$ be a stopping time less than the explosion time of $\widetilde{X}(\cdot)$, and, for an $\widetilde{\sr{F}}_{\tau}$-measurable random variable $U$ and a bounded measurable function $f$, define
\begin{align*}
   \varphi_{\tau,t}(U) = \dd{E}\left[ \left. f(|\widetilde{X}|(\tau+t)) \right| U \right], \qquad t \ge 0,
\end{align*}
then $\widetilde{Z}(\cdot)$ is a strong Markov process up to the explosion time if we show that
\begin{align}
\label{eq:c-ex 1}
  \dd{E}\left[ f(|\widetilde{X}|(\tau+t)) 1_{A} \right] = \dd{E}\left[\varphi_{\tau,t}(|\widetilde{X}|(\tau)) 1_{A} \right], \qquad A \in \widetilde{\sr{F}}_{\tau}, t > 0.
\end{align}
Since $\widetilde{X}(\cdot)$ is a strong Markov process and $|\widetilde{X}|(t) = |\widetilde{X}(t)|$, we have
\begin{align}
\label{eq:c-ex 2}
  \dd{E}\left[ f(|\widetilde{X}|(\tau+t)) 1_{A} \right] = \dd{E}\left[\varphi_{\tau,t}(\widetilde{X}(\tau)) 1_{A} \right], \qquad A \in \widetilde{\sr{F}}_{\tau}, t > 0,
\end{align}
We will show that \eq{c-ex 2} implies \eq{c-ex 1}. For this, note that \eq{tX-SIE} implies that
\begin{align}
\label{eq:t(-X)-SIE}
  - \widetilde{X}(t) & = - \widetilde{X}(0) - \int_{0}^{t} \widetilde{\sigma}(\widetilde{X}(s)) dW(s) - \int_{0}^{t} \widetilde{b}(\widetilde{X}(s)) ds \nonumber\\
  & = - \widetilde{X}(0) + \int_{0}^{t} \widetilde{\sigma}(-\widetilde{X}(s)) d(-W(s)) + \int_{0}^{t} \widetilde{b}(-\widetilde{X}(s)) ds, \qquad t \ge 0,
\end{align}
because $- \widetilde{b}(x) = - \sgn(x) b(|x|) = \sgn(-x) b(|-x|)$. Hence, $(-\widetilde{X})(\cdot) \equiv \{-\widetilde{X}(t); t \ge 0\}$ must be the weak solution of \eq{tX-SIE}. Hence, $\widetilde{X}(\cdot)$ and $(-\widetilde{X})(\cdot)$ are stochastically identical. This implies that the conditional distribution of $\{|\widetilde{X}|(t); t \ge \tau\}$ given $\widetilde{X}(\tau)$ is identical with that given $-\widetilde{X}(\tau)$. Namely, 
\begin{align*}
  \dd{E}[\varphi(\widetilde{X}(\tau)) 1_{A}] = \dd{E}[\varphi(-\widetilde{X}(\tau)) 1_{A}], \qquad A \in \widetilde{\sr{F}}_{\tau}.
\end{align*}
Hence, from \eq{c-ex 2},
\begin{align*}
 & \dd{E}\left[ f(|\widetilde{X}|(\tau+t)) 1_{A} \right] = \dd{E}\left[\varphi_{\tau,t}(\widetilde{X}(\tau)) 1_{A} \right] \nonumber\\
 & \quad = \dd{E}\left[\varphi_{\tau,t}(\widetilde{X}(\tau)) 1(\widetilde{X}(\tau) < 0) 1_{A} \right] + \dd{E}\left[\varphi_{\tau,t}(\widetilde{X}(\tau)) 1(\widetilde{X}(\tau) \ge 0)) 1_{A} \right] \nonumber\\
 & \quad = \dd{E}\left[\varphi_{\tau,t}(-\widetilde{X}(\tau)) 1(\widetilde{X}(\tau) < 0) 1_{A} \right] + \dd{E}\left[\varphi_{\tau,t}(\widetilde{X}(\tau)) 1(\widetilde{X}(\tau) \ge 0) 1_{A} \right] \nonumber\\
 & \quad = \dd{E}[\varphi(|\widetilde{X}|(\tau)) 1_{A}], \qquad A \in \widetilde{\sr{F}}_{\tau},
\end{align*}
because $\{\widetilde{X}(\tau) < 0\} \in \widetilde{\sr{F}}_{\tau}$. Thus, we have \eq{c-ex 1}, and the part (ii) is proved.

We finally prove (iii). For this, define $\widetilde{\beta}(x) = 2\widetilde{b}(x)/\widetilde{\sigma}^{2}(x)$ for $x \in \dd{R}$, and define
\begin{align*}
    \widetilde{\eta}(x) = \begin{cases}
  \int_{0}^{x} \exp\left(- \int_{0}^{y} \widetilde{\beta}(v) dv\right) dy, \quad & x \ge 0, \\
  - \int_{x}^{0} \exp\left(\int_{y}^{0} \widetilde{\beta}(v) dv\right) dy, \quad & x < 0.
\end{cases}
\end{align*}
Then, $\widetilde{\eta}(x) = \sgn(x) \eta(|x|)$ for $x \in R$ because
\begin{align*}
   \int_{x}^{0} \exp\left(\int_{y}^{0} \widetilde{\beta}(v) dv\right) dy = \int_{0}^{|x|} \exp\left( - \int_{0}^{y} \widetilde{\beta}(v) dv\right) dy, \qquad x < 0,
\end{align*}
and, obviously,
\begin{align*}
 & \widetilde{\eta}(y) - \widetilde{\eta}(x) = \int_{x}^{y} \exp\left(- \int_{0}^{w} \widetilde{\beta}(v) dv\right) dw > 0, \qquad x < y, \; x,y \in \dd{R}, \nonumber\\
 & \widetilde{\eta}'(x) = \exp\left(- \int_{0}^{x} \widetilde{\beta}(v) dv\right), \qquad \widetilde{\eta}''(x) = - \widetilde{\beta}(x) \exp\left(- \int_{0}^{x} \widetilde{\beta}(v) dv\right),
\end{align*}
where $\widetilde{\eta}''(x)$ is a Radon-Nikodym derivative. Using this fact, we apply the generalized version of Ito formula \eq{Ito1} to $\widetilde{X}(\cdot)$ for test function $\widetilde{\eta}$, we have, from \eq{tX-SIE},
\begin{align*}
  \widetilde{\eta}(\widetilde{X}(t)) & = \widetilde{\eta}(\widetilde{X}(s)) + \int_{s}^{t} \left(\widetilde{\eta}'(\widetilde{X}(u)) \widetilde{b}(\widetilde{X}(u)) + \frac 12 \widetilde{\eta}''(\widetilde{X}(u)) \widetilde{\sigma}^{2}(\widetilde{X}(u)) \right) du \nonumber\\
  & \qquad + \int_{s}^{t} \widetilde{\eta}'(\widetilde{X}(u)) \widetilde{\sigma}(\widetilde{X}(u)) \widetilde{W}(du) \nonumber\\
  & = \widetilde{\eta}(\widetilde{X}(s)) + \int_{s}^{t} \exp\left( - \int_{0}^{\widetilde{X}(u)} \widetilde{\beta}(v) dv \right) \widetilde{\sigma}(\widetilde{X}(u)) \widetilde{W}(du), \qquad 0 \le s < t,
\end{align*}
Hence, $\widetilde{M}(\cdot) \equiv \{\widetilde{\eta}(\widetilde{X}(t)); t \ge 0\}$ is a martingale with finite quadratic variations $[\widetilde{M}]_{t}$: 
\begin{align*}
  [\widetilde{M}]_{t} & = \int_{0}^{t} \exp\left( - 2 \int_{0}^{\widetilde{X}(u)} \widetilde{\beta}(v) dv \right) \widetilde{\sigma}^{2}(\widetilde{X}(u)) du, \qquad t \ge 0.
\end{align*}
Thus, let $c, d \in \dd{R}$ such that $c < d$, and let $\widetilde{\tau}_{x} = \inf\{t > 0; \widetilde{X}(t) = x\}$, then, by the optional sampling theorem, $\dd{E}_{x}[\widetilde{\eta}(\widetilde{X}(\widetilde{\tau}_{c} \wedge \widetilde{\tau}_{d}))] = \widetilde{\eta}(x)$ for $x \in (c,d)$, which yields
\begin{align*}
  \widetilde{\eta}(c) \dd{P}_{x}[ \widetilde{\tau}_{c} < \widetilde{\tau}_{d}] + \widetilde{\eta}(d) \dd{P}_{x}[ \widetilde{\tau}_{d} < \widetilde{\tau}_{c}] = \widetilde{\eta}(x).
\end{align*}
Since $\dd{P}_{x}[ \widetilde{\tau}_{c} < \widetilde{\tau}_{d}] + \dd{P}_{x}[ \widetilde{\tau}_{d} < \widetilde{\tau}_{c}] = 1$, we have, for $c < x < d$,
\begin{align}
\label{eq:tau-c-d}
  \dd{P}_{x}[ \widetilde{\tau}_{c} < \widetilde{\tau}_{d}] = \frac {\widetilde{\eta}(d) - \widetilde{\eta}(x)} {\widetilde{\eta}(d) - \widetilde{\eta}(c)}, \qquad \dd{P}_{x}[ \widetilde{\tau}_{d} < \widetilde{\tau}_{c}] = \frac {\widetilde{\eta}(x) - \widetilde{\eta}(c)} {\widetilde{\eta}(d) - \widetilde{\eta}(c)}.
\end{align}
Since $\lim_{c \downarrow -\infty} \widetilde{\tau}_{c} = \lim_{d \uparrow \infty} \widetilde{\tau}_{d} = \infty$ almost surely and $\widetilde{\eta}(x)$ is strictly increasing, \eq{tau-c-d} implies that
\begin{align*}
  \dd{P}_{x}[ \widetilde{\tau}_{c} < \infty] = \lim_{d \uparrow \infty} \dd{P}_{x}[ \widetilde{\tau}_{c} < \widetilde{\tau}_{d}], \qquad \dd{P}_{x}[ \widetilde{\tau}_{d} < \infty] = \lim_{c \downarrow -\infty} \dd{P}_{x}[ \widetilde{\tau}_{d} < \widetilde{\tau}_{c}].
\end{align*}
Hence, $\dd{P}_{x}[ \widetilde{\tau}_{y} < \infty] = 1$ for any $x,y \in \dd{R}$ if and only if \eq{eta-infty} holds because $\widetilde{\eta}(x) = \sgn(x) \eta(|x|)$. Hence, $\widetilde{X}(\cdot)$ is recurrent if and only if \eq{eta-infty}. Because $\widetilde{Z}(\cdot) = |\widetilde{X}|(\cdot)$, the same property holds for $\widetilde{Z}(\cdot)$. This proves (iii).
\end{proof}

\section{Stationary measure and distribution}
\label{sec:stationary}
\setnewcounter

By (iii) of \lem{Z-X-existence}, under the condition \eq{eta-infty}, $Z(\cdot)$ has a unique stationary measure $\nu$, that is, for any bounded measurable function $f$ from $\dd{R}_{+}$ to $\dd{R}$,
\begin{align*}
  \int_{0}^{\infty} f(y) \nu(dy) = \int_{0}^{\infty} \dd{E}_{y}[f(Z(t))] \nu(dy).
\end{align*}
We show that this $\nu$ is obtained in a closed form.

\begin{theorem}
\label{thr:stationary}
Assume \cond{Z-existence} and \eq{eta-infty}. Then, the solution $Z(\cdot)$ of \eq{Z-SIE} obtained in \lem{Z-X-existence} has a stationary measure $\nu$ which is absolutely continuous with respect to Lebesgue measure on $[0,\infty)$ and has the following density function $h(x)$ up to a constant multiplication.
\begin{align}
\label{eq:stationary 1}
 & h(x) = \frac {1}{\sigma^{2}(x)} \exp\left(\int_{0}^{x} \beta(u) du\right), \qquad x \in \dd{R}_{+}, \qquad x \ge 0,
\end{align}
where recall $\beta(x) = 2 b(x)/\sigma^{2}(x)$ of \eq{beta}. In particular, if
\begin{align}
\label{eq:C 1}
 & C_{[0,\infty)} \equiv \int_{0}^{\infty} \frac {1}{\sigma^{2}(y)} \exp\left(\int_{0}^{y} \beta(u) du\right) dy < \infty,
\end{align}
then there is the unique stationary distribution which has the density $C_{[0,\infty)}^{-1} h(x)$ for $x \ge 0$. Furthermore, there is a probability law for $Z(\cdot)$ to be a stationary process, and, under this probability law,
\begin{align}
\label{eq:EY}
  \dd{E}[Y(1)] = \frac 12 \dd{E}[L_{0}(1)] = \frac 12 C_{[0,\infty)}^{-1}.
\end{align}

 \end{theorem}
 
\begin{remark}
\label{rem:stationary}
As discussed in \rem{Z-uniqueness}, the weak uniqueness of $Z(\cdot)$ is not known, but this theorem holds true as long as $Z(\cdot)$ has a stationary measure. Hence, the stationary measure and distribution of $Z(\cdot)$ uniquely exist as long as they exist.
\end{remark}

\begin{proof}
We apply the generalized Ito formulas \eq{G-Ito-X} with a convex test function to the $Z(\cdot)$, then we have \eq{G-Ito-Z}. 
Define a convex test $f$ for fixed $x \ge 0$ as
\begin{align}
\label{eq:SIE-test1}
  f(u) = (x-u)^{+} 1(u < x). \qquad u\in \dd{R}.
\end{align}
Since $f'(x+) = 0$ and $f'(x-) = -1$, it follows from \eq{Mf-1} that
\begin{align*}
  \mu_{f}(\{x\}) = \lim_{\varepsilon \downarrow 0} f'((x+\varepsilon) -) - f'(x-) = f'(x+) - f'(x-) = 1.
\end{align*}
On the other hand, $f'(u) = - 1(u < x)$ and $f''(u-) = 0$ for $u < x$. Hence,
\begin{align*}
 & \int_{0}^{\infty} L_{u}(t) \mu_{f}(du) = L_{x}(t),
\end{align*}
and the generalized Ito formula \eq{G-Ito-Z} for the $Z(\cdot)$ yields
\begin{align}
\label{eq:G-Ito1}
  f(Z(t)) & = f(Z(0)) - \int_{0}^{t} 1(0 \le Z(s) < x) \sigma(Z(s)) dW(s)  \nonumber\\
 & \quad - \int_{0}^{t} b(Z(s)) 1(0 \le Z(s) < x) ds - Y(t) + \frac 12 L_{x}(t), \qquad t \ge 0.
\end{align}
Let $x=0$ for this equation, then we have
\begin{align}
\label{eq:Y-L0}
  f(Z(t)) = f(Z(0)) - Y(t) + \frac 12 L_{0}(t), \qquad t \ge 0.
\end{align}

Applying \eq{Y-L0} to \eq{G-Ito1}, and take the conditional expectation of each term in it for $t=1$ given $Z(0)=y$, we have
\begin{align}
\label{eq:G-Ito2}
    \int_{0}^{1}  \dd{E}_{y}\left[b(Z(s)) 1(0 \le Z(s) < x)\right] ds = - \frac 12 \dd{E}_{y}[L_{0}(1)] + \frac 12 \dd{E}_{y}[L_{x}(1)], \qquad t \ge 0.
\end{align}
We integrate both sides of this equation by the invariant measure $\nu$ of $Z(\cdot)$. For this, we use \eq{L-time1} for $t=1$ and $p(u) = 1(0 \le u \le x)/\sigma^{2}(u)$, which yields
\begin{align}
\label{eq:L-time5}
  \int_{0}^{x} \dd{E}_{y}[L_{u}(1)] \frac 1{\sigma^{2}(u)} du = \int_{0}^{1} \dd{E}_{y}[p(Z(s)) \sigma^{2}(Z(s))] ds = \int_{0}^{1} \dd{P}_{y}[0 \le Z(s) \le x] ds.
\end{align}
Define
\begin{align}
\label{eq:m-nu}
  m_{\nu}[L_{x}(1)] = \int_{0}^{\infty} \dd{E}_{y}[L_{x}(1)] \nu(dy), \qquad x \ge 0.
\end{align}
Then, integrating \eq{L-time5} by the stationary measure $\nu$, we have
\begin{align}
\label{eq:L-time6}
  \int_{0}^{x} m_{\nu}[L_{u}(1)] \frac 1{\sigma^{2}(u)} du =  \int_{0}^{1} \int_{0}^{\infty} \dd{P}_{y}[0 \le Z(s) \le x] \nu(dy) ds = \nu([0,x]) < \infty.
\end{align}
Since $1/\sigma^{2}(u)$ is locally integrable by \eq{local integrability 2} and $L_{u}(1)$ is continuous in $u$, $m_{\nu}[L_{u}(1)]$ must be finite for any $u \ge 0$.
Furthermore, since $f$ is bounded function on $[0,\infty)$ with bounded support $[0,x]$ for each fixed $x \ge 0$, integrating \eq{G-Ito2} by $\nu$, we have
\begin{align}
\label{eq:G-Ito1*}
 & \int_{0}^{x} b(y) \nu(dy) + \frac 12 \left(m_{\nu}[L_{0}(1)] - m_{\nu}[L_{x}(1)]\right) = 0.
\end{align}

We next relate $\nu$ to $\dd{E}[L_{x}(1)]$. For this, we again use \eq{L-time1} for $t=1$ similar to \eq{L-time6} but for $p(u) = \beta(u) 1(u \le x)$. Then, we have, by \eq{G-Ito1*}, for $x \ge 0$,
\begin{align}
\label{eq:IE-Lx 2}
  \int_{0}^{x} \beta(u) m_{\nu}[L_{u}(1)] du = \int_{0}^{1} \left(\int_{0}^{x} 2 b(y) \nu(dy)\right) ds = m_{\nu}[L_{x}(1)] - m_{\nu}[L_{0}(1)].
\end{align}
Hence, $m_{\nu}[L_{x}(1)]$ has the derivative $v(x) \equiv \beta(x) m_{\nu}[L_{x}(1)]$ almost surely with respect to Lebesgue measure on $[0,\infty)$, and therefore \eq{IE-Lx 2} yields
\begin{align}
\label{eq:DE 1}
   \frac d{dx} \left( \exp\left(-\int_{0}^{x} \beta(y) dy\right) m_{\nu}[L_{x}(1)] \right) = 0, \qquad x \ge 0.
\end{align}
This concludes
\begin{align}
\label{eq:Lx}
  m_{\nu}[L_{x}(1)] = m_{\nu}[L_{0}(1)] \exp\left( \int_{0}^{x} \beta(y) dy \right), \qquad x \ge 0.
\end{align}

Finally, plugging \eq{Lx} into \eq{L-time6} with $p(x) = 1/\sigma^{2}(x)$, we have
\begin{align*}
  \nu([0,x]) = m_{\nu}[L_{0}(1)] \int_{0}^{x} \frac {1} {\sigma^{2}(y)} \exp\left( \int_{0}^{y} \beta(u) du \right) dy, \qquad x \ge 0.
\end{align*}
This proves \eq{stationary 1}. If \eq{C 1} holds, then $\nu$ is a finite measure, and can be normalized by multiplying $(C_{[0,\infty)} m_{\nu}[L_{0}(1)])^{-1}$. Thus, the stationary distribution has the density $C_{[0,\infty)}^{-1} h(x)$. Taking this distribution for $Z(0)$, $Z(\cdot)$ is stationary process, and \eq{EY} is immediate from \eq{Y-L0}.
\end{proof}

\thr{stationary} generalizes Theorem 3.2 of \citet{Miya2024b}, which assumes that $\sigma(x)$ and $\sigma(x)$ are simple functions. We discuss the literature about the reflecting diffusions in \sectn{concluding}.

\section{Extension to the two-sides reflected diffusion}
\label{sec:extension}
\setnewcounter

In this section, we consider a nonnegative diffusion process reflected not only at the origin but also at $a > 0$. We define this two-sides reflected diffusion as the solution of $U(\cdot) \equiv \{U(t); t \ge 0\}$ of the stochastic integral equation below.
\begin{align}
\label{eq:U-SIE}
  U(t) = U(0) & + \int_{0}^{t} \sigma(U(s)) dW(s) + \int_{0}^{t} b(U(s)) ds + Y_{0}(t) - Y_{a}(t) \in [0,a], \; t \ge 0,
\end{align}
where $Y_{x}(\cdot)\equiv \{Y_{x}(t); t \ge 0\}$ is a non-deceasing process for $x=0, a$ satisfying that
\begin{align*}
  \int_{0}^{t} 1(Z(s) > 0) dY_{0}(s) = 0, \qquad \int_{0}^{t} 1(Z(s) < a) dY_{a}(s) = 0, \qquad t \ge 0.
\end{align*}

For constructing the weak solution of \eq{U-SIE}, we apply the same idea which was used for constructing the solution $Z(\cdot)$ of SIE \eq{Z-SIE}.  Since \eq{U-SIE} does not depend on $\sigma(x)$ and $b(x)$ for $x \in \dd{R} \setminus [0,a]$, we further weaken \cond{Z-existence} to
\begin{condition}
\label{cond:U-existence}
The functions $b(x)$ and $\sigma(x)$ defined on $[0,a]$ satisfy
\begin{align}
\label{eq:sigma positive 3}
 & \sigma(x) > 0, \qquad \forall x \in [0,a],\\
\label{eq:local integrability 4}
 & \int_{x_{1}}^{x_{2}} \frac 1{\sigma^{2}(y)} dy < \infty, \qquad \forall (x_{1},x_{2}) \in [0,a] \mbox{ satisfying } x_{1} < x_{2},\\
\label{eq:local integrability 5}
 & \int_{0}^{x+\varepsilon} |\beta(y)| dy < \infty, \qquad \forall x \in [0,a), \exists \varepsilon \in (0,a - x].
\end{align}
\end{condition}

Similar to \lem{Z-X-existence}, we have
\begin{lemma}
\label{lem:U-X-existence}
Under \cond{U-existence}, (i) the stochastic integral equation \eq{U-SIE} has a weak solution $({U}(\cdot),{Y}_{0}(\cdot),{Y}_{a}(\cdot),{W}(\cdot))$ up to the explosion time of $U(\cdot)$, and (ii) this ${U}(\cdot)$ is a strong Markov process up to its explosion time, where the same notations are used for the weak solution for simplicity. Furthermore, (iii) the ${U}(\cdot)$ is Harris recurrent.
\end{lemma}

Based on this lemma, we derive the stationary distribution of $U(\cdot)$ similarly to \thr{stationary}.

\begin{theorem}
\label{thr:stationary-a}
Assume \cond{U-existence}. Then, the weak solution $U(\cdot)$ of \eq{U-SIE} obtained in \lem{U-X-existence} has a stationary measure $\nu$ which is absolutely continuous with respect to Lebesgue measure on $[0,\infty)$ and has the following density function.
\begin{align}
\label{eq:stationary-a}
 & h(x) = \frac {1}{\sigma^{2}(x)} \exp\left(\int_{0}^{x} \beta(u) du\right), \qquad x \in [0,a],
\end{align}
where recall $\beta(x) = 2 b(x)/\sigma^{2}(x)$ of \eq{beta}. Particularly, if
\begin{align}
\label{eq:C-a 1}
 & C_{[0,a]} \equiv \int_{0}^{a} \frac {1}{\sigma^{2}(y)} \exp\left(\int_{0}^{y} \beta(u) du\right) dy < \infty,
\end{align}
then $U(\cdot)$ has the stationary distribution whose density is $C_{[0,a]}^{-1} h(x)$ for $x \in [0,a]$. Furthermore, under the probability law for $U(\cdot)$ to be stationary,
\begin{align}
\label{eq:EY0-Ya}
  \dd{E}[Y_{0}(1)] = \frac 12 C_{[0,a]}^{-1}, \qquad \dd{E}[Y_{a}(1)] = \frac 12 C_{[0,a]}^{-1} \exp\left(\int_{0}^{a} \beta(u) du\right).
\end{align}
\end{theorem}

\begin{remark}
\label{rem:uniqueness-a}
Similarly to \thr{stationary}, the stationary measure is uniquely given by \eq{stationary-a} up to a multiplication constant for any weak solution $U(\cdot)$ of \eq{U-SIE} if it has a stationary measure. $C_{[0,a]}$ of \eq{C-a 1} may be infinite because $\int_{0}^{a} |\beta(y)| dy = \infty$ is possible under \eq{local integrability 5}.
\end{remark}

\begin{remark}
\label{rem:stationary-a}
If the diffusion coefficients are sufficiently smooth and if the diffusion process $X(\cdot)$ is positive recurrent on $(a,b)$, then the same stationary distribution is obtained under much stronger assumptions in Theorem 12.2 of \cite{BhatWaym2009}. The background of this type of results is discussed in \sectn{concluding}. 
\end{remark}

 In the rest of this section, we prove these lemma and theorem.

\begin{proof}[Proof of \lem{U-X-existence}]
(i) We use a similar idea for the proof of (i) of \lem{Z-X-existence}. Namely, we extend the $\sigma(x)$ and $b(x)$ on $[0,a]$ to diffusion coefficients on $\dd{R}$ so that they symmetrically turn up at $x=a$. For this extension, we introduce the following functions.
\begin{align*}
&  g_{a}(x) = (x-a) \sgn(a-x) + a, \quad g_{0}(x) = (|x-a| -a) 1(|a-x| > a), \quad x \in \dd{R},
\end{align*}
and define $f$ as
\begin{align*}
  g(x) = g_{a}(x) + g_{0}(x), \qquad x \in \dd{R}.
\end{align*}
See \fig{test-function} for these functions. 
\begin{figure}[h] %  figure placement: here, top, bottom, or page
   \centering
   \includegraphics[width =10cm]{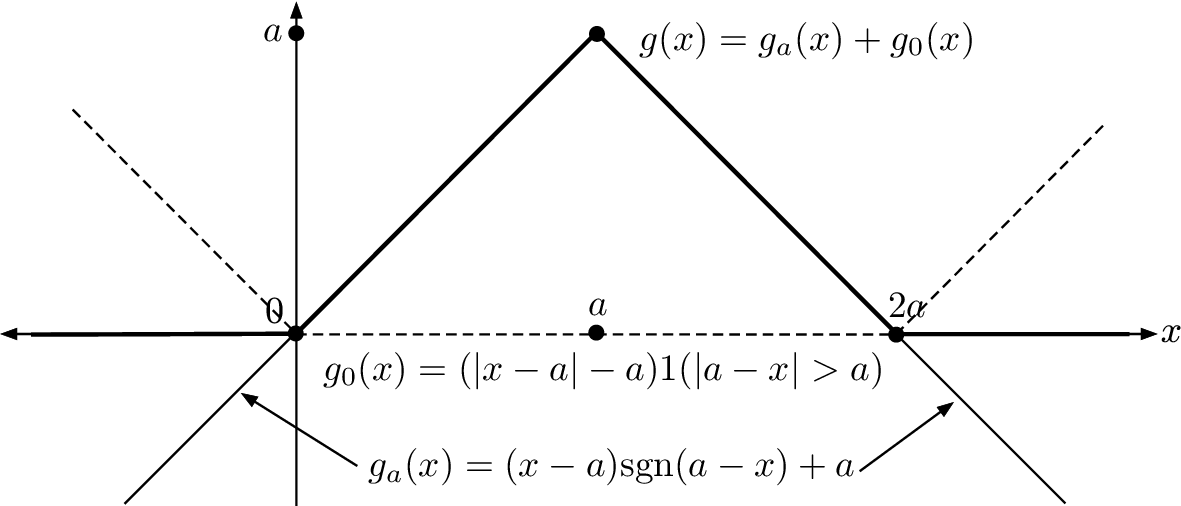} %\vspace{-1ex}
   \caption{Test function $g$}
   \label{fig:test-function}
\end{figure}

We define the extensions of $\sigma(x)$ and $b(x)$ on $[0,a]$ as
\begin{align*}
%\label{eq:ora sigma-b}
  \widehat{\sigma}(x) = \begin{cases}
  1, & x \in \dd{R} \setminus [0,2a], \\
  \sigma(g(x)), & x \in [0,2a],
\end{cases}, \quad \widehat{b}(x) = \begin{cases}
  -1, & x \in \dd{R} \setminus [0,2a],\\
  \sgn(a - x) b(g(x)), &  x \in [0,2a].
\end{cases}
\end{align*}
Obviously, \cond{X-existence} is satisfied for $\widehat{\sigma}(x)$ and $\widehat{b}(x)$ under \cond{U-existence}. Hence, by \lem{X-existence}, $\widehat{X}(\cdot) = \{\widehat{X}(t); t \ge 0\}$ can be defined as the weak solution the following stochastic integral equation. 
\begin{align}
\label{eq:hX-SIE}
  \widehat{X}(t) = \widehat{X}(0) + \int_{0}^{t} \widehat{\sigma}(\widehat{X}(s)) dW(s) + \int_{0}^{t} \widehat{b}(\widehat{X}(s)) ds, \qquad t \ge 0.
\end{align}
Clearly, $g_{0}$ and $g_{a}$ are concave and convex functions, respectively. Hence, the test function $g$ is the difference of convex functions, so it can be used as a test function for the generalized Ito formula. For computing it, we note that
\begin{align*}
 & g'(x) = \begin{cases}
 1, \quad & x \in (0,a), \\
 -1, &  x \in (a,2a),\\
 0, & \mbox{ otherwise},
\end{cases},\\
 & \mu_{g_{0}}(B) = 1(0 \in B) - 1(2a \in B), \\ & \mu_{-g_{a}}(B) = 2 \times 1(a \in B), \qquad B \in \sr{B}(\dd{R}).
\end{align*}
Since $g(x) = g_{0}(x) - (-g_{a}(x))$, for the weak solution $X(\cdot)$ of \eq{X-SIE}, it follows from the generalized Ito formula that
\begin{align}
\label{eq:ora(X) 1}
 & g(\widehat{X}(t)) - g(\widehat{X}(0)) = \int_{0}^{t} g'(\widehat{X}(s)) d\widehat{X}(s) + \frac 12 \int_{0}^{\infty} \widehat{L}_{x}(t) \mu_{g_{0}}(dx) - \frac 12 \int_{0}^{\infty} \widehat{L}_{x}(t) \mu_{-g_{a}}(dx) \nonumber\\
  & \quad = \int_{0}^{t} g'(\widehat{X}(s)) \left(\widehat{\sigma}(\widehat{X}(s)) dW(s) + \widehat{b}(\widehat{X}(s)) ds \right) + \frac 12 \left( \widehat{L}_{0}(t) + \widehat{L}_{2a}(t)\right) - \widehat{L}_{a}(t) \nonumber\\
  & \quad = \int_{0}^{t} \sigma(g(\widehat{X}(s))) [1(\widehat{X}(s) \in (0,a)) - 1(\widehat{X}(s) \in (a,2a))] dW(s) + \int_{0}^{t} b(g(\widehat{X}(s))) ds \nonumber\\
  & \qquad  + \frac 12 \left( \widehat{L}_{0}(t) + \widehat{L}_{2a}(t)\right) - \widehat{L}_{a}(t).
\end{align}
Define $\widehat{W}(\cdot) \equiv \{\widehat{W}(t); t \ge 0$ as
\begin{align*}
  \widehat{W}(t) = \sgn(a - \widehat{X}(t)) W(t), \qquad t \ge 0,
\end{align*}
then $\widehat{W}(\cdot)$ stochastically identical with $W(\cdot)$, and, similar to \eq{substitution-f},
\begin{align*}
  \int_{0}^{t} \sigma(|\widehat{X}|(s)) \sgn(\widehat{X}(s)) dW(s) = \int_{0}^{t}  \sigma(|\widehat{X}|(s)) d\widehat{W}(s).
\end{align*}
Furthermore, $\widehat{L}_{0}(t) = \widehat{L}_{2a}(t)$ because
\begin{align*}
  \int_{-\infty}^{+\infty} \widehat{L}_{x}(t) p(x) dx = \int_{0}^{t} p(\widehat{X}(s)) \sigma^{2}(g_{a}(\widehat{X}(s))) ds
\end{align*}
and $g_{a}(\varepsilon) = - g_{a}(2a-\vc{\varepsilon})$ for $\varepsilon [0,a)$. Hence, \eq{ora(X) 1} can be written as
\begin{align}
\label{eq:ora(X) 2}
  g(\widehat{X}(t)) & = g(\widehat{X}(0))  \nonumber\\
  & \qquad + \int_{0}^{t} \sigma(g(\widehat{X}(s))) d\widehat{W}(s) + \int_{0}^{t} b(g(\widehat{X}(s))) ds + \widehat{L}_{0}(t) - \widehat{L}_{a}(t).
\end{align}
Consequently, let
\begin{align*}
  U(t) = g(\widehat{X}(t)), \qquad Y_{0}(t) = \widehat{L}_{0}(t), \qquad Y_{a}(t) = \widehat{L}_{a}(t),
\end{align*}
then, we have the weak solution $(U(\cdot), Y_{0}(\cdot), Y_{a}(\cdot), \widehat{W}(\cdot))$ of SIE \eq{U-SIE}. 

(ii) The strong Markov property can be proved in the exactly same ways as in the proof of \lem{Z-X-existence}, so we omit its proof.

(iii) The recurrence condition of $U(\cdot)$ is proved similarly to the proof of \lem{Z-X-existence}, but some changes are needed because $\widetilde{\sigma}(x)$ and $\widetilde{b}(x)$ are changed to $\widehat{\sigma}(x)$ and $\widehat{b}(x)$. So, we outline the proof. Define $\widehat{\beta}(x) = 2\widehat{b}(x)/\widehat{\sigma}^{2}(x)$ for $x \in \dd{R}$, and define
\begin{align*}
    \widehat{\eta}(x) = \begin{cases}
 \int_{0}^{x} \exp\left(- \int_{0}^{y} \widehat{\beta}(v) dv\right) dy, \quad & x \ge 0,\\
 - \int_{0}^{-x} \exp\left(- \int_{0}^{y} \widehat{\beta}(v) dv\right) dy, \quad & x < 0.
\end{cases}
\end{align*}
Obviously, 
\begin{align}
\label{eq:hbeta}
  \widehat{\beta}(x) = \begin{cases}
  -2, & x \in \dd{R} \setminus [0,a],\\
  \sgn(a-x) \beta(g(x)), \quad & x \in [0,a],   
\end{cases}
\end{align}
$\widehat{\eta}(s)$ is increasing in $x \in \dd{R}$, and $\widehat{\eta}''(x) + \widehat{\beta}(x)  \widehat{\eta}'(x) = 0$ almost surely.
Using these facts, we apply the generalized version of Ito formula \eq{Ito1} to $\widehat{X}(\cdot)$ for test function $\widehat{\eta}$, we have, from \eq{hX-SIE},
\begin{align*}
  \widehat{\eta}(\widehat{X}(t)) & = \widehat{\eta}(\widehat{X}(s)) + \int_{s}^{t} \exp\left( - \int_{0}^{\widehat{X}(u)} \widehat{\beta}(v) dv \right) \widehat{\sigma}(\widehat{X}(u)) \widehat{W}(du), \qquad 0 \le s < t,
\end{align*}
Hence, $\widehat{M}(\cdot) \equiv \{\widehat{\eta}(\widehat{X}(t)); t \ge 0\}$ is a martingale with finite quadratic variations. Thus, let $c, d \in \dd{R}$ such that $c < d$, and let $\widehat{\tau}_{x} = \inf\{t> 0; \widehat{X}(t) = x\}$, then, similarly to the proof of \lem{Z-X-existence}, we have, for $c < x < d$,
\begin{align}
\label{eq:tau-h-c-d}
  \dd{P}_{x}[ \widehat{\tau}_{c} < \widehat{\tau}_{d}] = \frac {\widehat{\eta}(d) - \widehat{\eta}(x)} {\widehat{\eta}(d) - \widehat{\eta}(c)}, \qquad \dd{P}_{x}[ \widehat{\tau}_{d} < \widehat{\tau}_{c}] = \frac {\widehat{\eta}(x) - \widehat{\eta}(c)} {\widehat{\eta}(d) - \widehat{\eta}(c)}.
\end{align}
Since \eq{hbeta} implies that $\lim_{d \uparrow \infty} \widehat{\eta}(x) = \lim_{c \downarrow - \infty} \widehat{\eta}(x) = \infty$, \eq{tau-c-d} implies that
\begin{align*}
  \dd{P}_{x}[ \widehat{\tau}_{c} < \infty] = \lim_{d \uparrow \infty} \dd{P}_{x}[ \widehat{\tau}_{c} < \widehat{\tau}_{d}] = 1, \qquad \dd{P}_{x}[ \widehat{\tau}_{d} < \infty] = \lim_{c \downarrow -\infty} \dd{P}_{x}[ \widehat{\tau}_{d} < \widehat{\tau}_{c}] = 1.
\end{align*}
Hence, $\dd{P}_{x}[ \widehat{\tau}_{y} < \infty] = 1$ for any $x,y \in [0,a]$, and therefore $\widehat{X}(\cdot)$ is recurrent. Because $\widehat{Z}(t) = g (\widehat{X}(t))$, the same property holds for $\widehat{Z}(\cdot)$. This proves (iii).
\end{proof}

\begin{proof}[Proof of \thr{stationary-a}]
This theorem is proved essentially in the same way as the proof of \thr{stationary}, applying the generalized Ito formulas \eq{G-Ito-X} with the same test function $f$ whose domain is restricted on $[0,a]$. Namely, for $x \in [0,a)$,
\begin{align*}
  f(u) = (x-u)^{+} 1(u < x). \qquad u \in [0,a].
\end{align*}
Then, we have the same formula as \eq{G-Ito2} in which $Z(\cdot)$ is replaced by $U(\cdot)$ because the integration by $Y_{a}(t)$ vanishes. Hence, the same arguments as in the proof of \thr{stationary} go through, and we have \eq{stationary-a}. It remains to prove \eq{EY0-Ya}. For this, we apply Ito formula to \eq{U-SIE} for a continuously twice differentiable function $f$, and take its expectation under the probability measure for $U(\cdot)$ to be stationary. Then, we have
\begin{align*}
  \int_{0}^{1} \dd{E} \left[f'(U(s)) b(U(s)) + \frac 12 f''(U(s)) \sigma^{2}(U(s)) \right]ds + f'(0) \dd{E}[Y_{0}(1)] - f'(a) \dd{E}[Y_{a}(1)] = 0.
\end{align*}
Since $U(s)$ has the density function $C_{[0,a]}^{-1} h(x)$ on $[0,x]$, this yields
\begin{align*}
  \frac 1{2C} \int_{0}^{a} \left(f'(x) \beta(x) + f''(x)\right) \exp\left( \int_{0}^{x} \beta(u) du\right) dx + f'(0) \dd{E}[Y_{0}(1)] - f'(a) \dd{E}[Y_{a}(1)] = 0.
\end{align*}
Since $\left(f'(x) \beta(x) + f''(x)\right) \exp\left( \int_{0}^{x} \beta(u) du\right) = \left(f'(x) \exp\left( \int_{0}^{x} \beta(u) du\right)\right)'$, this can be written as
\begin{align*}
  f'(0) \left(\dd{E}[Y_{0}(1)] - \frac 1{2C_{[0,a]}}\right) + f'(a) \left( \frac 1{2C_{[0,a]}} \exp\left( \int_{0}^{a} \beta(u) du\right) - \dd{E}[Y_{a}(1)]\right) = 0.
\end{align*}
Choosing $f$ such that $f'(0) = 0$ and $f'(a) \not= 0$, we have \eq{EY0-Ya}.
\end{proof}

\section{Concluding remarks}
\label{sec:concluding}
\setnewcounter

In this section, we discuss about the literature on one-dimensional reflected diffusions and their stationary distributions.

We first note that the reflected diffusion $Z(\cdot)$ and its stationary distribution have not been well studied when the diffusion coefficients are discontinuous except for \citet{Miya2024b}. This is because they may be out of interest as an exceptional case. However, as shown in \cite{Miya2025}, the discontinuous case may be also interesting in application. This is contrasted with the fact that diffusions including reflection have been widely studied when their diffusion coefficients are smooth, typically, Lipschitz continuous. For example, if the diffusion coefficients are Lipschitz continuous and bounded by $K\sqrt{(1+|x|^{2})}$ with a constant $K$, then \citet{Tana1979} shows that the strong solution of SIE \eq{Z-SIE} and its multi-dimensional versions uniquely exist even when the diffusion coefficients depend on time. As for their detailed analysis, one-dimensional reflected diffusions have been less studied even when the diffusion coefficients are continuous. For example, \citet{ReedWardZhan2013} study a special case of $Z(\cdot)$ under the name, the generalized drift Skorokhod problem solution, which assumes that $\sigma(x)$ is a constant and $b(x)$ is Lipschitz continuous, but its stationary distribution is not studied.

So far, in the literature, a different formulation is employed to study the one-dimensional reflected diffusions and their stationary distributions, particularly when the diffusion coefficients are discontinuous. As mentioned in \sectn{introduction}, this formulation firstly studies the diffusion process $X(\cdot)$ under \cond{X-existence}, and considers its stationary measure. Then, the restriction of this stationary measure on the half line or a closed (or open) interval is considered to be the stationary measure of the one or two sides reflected diffusion, respectively. This kind of arguments also are employed even for the case that the diffusion coefficients of $X(\cdot)$ are sufficiently smooth. We detail both of them below.

As for the general case which allows the diffusion coefficients to be discontinuous, we refer to Chapter 23 of \cite{Kall2001}, which only considers the case that the drift is null. Extend the scale function $\eta(x)$ of \eq{scale-f} to all $x \in \dd{R}$ by
\begin{align*}
  \eta(x) = - \int_{x}^{0} \exp\left(\int_{y}^{0} \beta(v) dv\right) dy, \qquad x < 0,
\end{align*}
where recall $\beta(x) = 2 b(x)/\sigma^{2}(x)$ for $x \in \dd{R}$. We call this $\eta$ the scale function of $X(\cdot)$, which is different from $\widetilde{\eta}$ in the proof of \lem{Z-X-existence}. Then, in the null drift case, $\eta(x) = x$, and therefore one can see that $X(\cdot)$ is recurrent by the same arguments as in the proof of \lem{Z-X-existence} if the explosion time $S$ of $X(\cdot)$ is infinite. Thus, if $\dd{P}[S = \infty] = 1$, then the $X(\cdot)$ with null drift always has a stationary measure. Restricting this measure on the half line and closed and open intervals, the main results of \cite{Kall2001} are obtained in Theorem 23.15 and Lemma 23.19, in which the condition $S=\infty$ is not specified, but must be. Under this condition, we can see that the null drift condition is  not essential because such a diffusion is obtained as $(\eta \circ X)(\cdot) \equiv \{\eta(X(t)); t \ge 0\}$, so the stationary measure of $(\eta \circ X)(\cdot)$ is transformed to that of $X(\cdot)$. For this transformation, we need to know the stationary measure of the null drift case, whose density is noticed to be $1/\sigma^{2}(x)$ without proof in \cite{Kall2001}. Since $\dd{P}[S = \infty] = 1$ is equivalent to
\begin{align}
\label{eq:eta-+infty}
  \lim_{x \downarrow -\infty} \eta(x) = - \infty, \qquad \lim_{x \uparrow \infty} \eta(x) = \infty,
\end{align}
we have the following lemma. We also can directly proof this lemma using similar arguments to the proof of \thr{stationary}.

\begin{lemma}\rm
\label{lem:X-stationary}
Assume \cond{X-existence}. Then, the unique weak solution $X(\cdot)$ of \eq{X-SIE} is recurrent if and only if \eq{eta-+infty} holds. Furthermore, in this recurrent case, $X(\cdot)$ has a stationary measure $\nu$ which is absolutely continuous with respect to Lebesgue measure on $[0,\infty)$ and has the following density function $h(x)$ up to a constant multiplication.
\begin{align}
\label{eq:X-stationary 1}
 & h(x) = \begin{cases}
  \frac {1}{\sigma^{2}(x)} \exp\left(\int_{0}^{x} \beta(u) du\right), \quad & x \ge 0, \\
  \frac {1}{\sigma^{2}(x)} \exp\left(\int_{-x}^{0} \beta(-u) du\right), \quad & x < 0. 
\end{cases}
\end{align}
In particular, if $C_{\dd{R}} \equiv \int_{-\infty}^{\infty} h(x) dx < \infty$ then there is the unique stationary distribution which has the density $C_{\dd{R}}^{-1} h(x)$ for $x \ge 0$.
\end{lemma}

We next consider the case that the diffusion coefficients are sufficiently smooth, typically, Lipschitz continuous. In this case, we refer to Chapter V of \cite{BhatWaym2009}, which defines the diffusion $X(\cdot)$ as a Markov process assuming the smooth condition (2.1)' of \cite{BhatWaym2009}, which is close to Lipschitz continuity. Restricting the diffusion $X(\cdot)$ on open interval $(c,d)$, the stationary distribution of this restricted process is derived in Theorem 12.2 of \cite{BhatWaym2009} assuming the existence of the stationary distribution and that the diffusion coefficients are twice continuously differentiable, which are much stronger than the assumptions in \thr{stationary-a} of this paper. This stronger continuity is required to have Kolmogorov forward equation, which is used to derive the stationary distribution. Without loss of generality, these $c,d$ in Theorem 12.2 can be replaced by $0,a$. Then, this stationary distribution is identical with the one in \thr{stationary-a} for $U(\cdot)$ on $[0,a]$. 

Chapter V of \cite{BhatWaym2009} also consider a reflected diffusion. For example, it is defined as a Markov process satisfying certain conditions in Definition 6.1 of \cite{BhatWaym2009}, in which the reflecting point $0$ is replaced by $a \in \dd{R}$, but this is not essential. If this Markov process is uniquely determined by those conditions, then we can see that the reflected diffusion of this definition indeed exists, and is identical with our $Z(\cdot)$. However, this uniqueness is not argued in \cite{BhatWaym2009}. Furthermore, it does not discuss about its stationary distribution.

Thus, the reflected diffusions themselves are not really defined in \cite{BhatWaym2009,Kall2001}, and it is unclear for what reflected diffusions their results are obtained. However, it is notable that their results are fully compatible with those obtained from the reflected diffusion $Z(\cdot)$. From this, we can see that there is no difference in the stationary distributions of the $X(\cdot)$ on $(c,d)$ and of the two sides reflected diffusion on $[c,d]$ as long as those stationary distribution exist. Namely, reflections in one dimensional diffusion do not influence the stationary distributions. Of course, this can not be expected for the multi-dimensional case.

We have only highlighted theoretical results for the stationary distributions, but there are their heuristic derivations (e.g., see \cite{BrowWhit1995}). They may be checked by Theorems \thrt{stationary} and \thrt{stationary-a}.

\appendix

\section*{Appendix}
\setnewcounter
\setcounter{section}{1}

%\section{Local time and generalized Ito formula}
%\label{app:framework}
\setnewcounter

\subsection{Local time}
\label{app:local-time}

We briefly discuss about local time for a semi-martingale. Here, a real-valued continuous-time stochastic process $X(\cdot) = \{X(t); t \ge 0\}$ is called a semi-martingale if it is decomposed as 
\begin{align*}
  X(t) = X(0) + M(t) + A(t), \qquad t \ge 0,
\end{align*}
for a martingale $M(\cdot) = \{M(t); t \ge 0\}$ and a process $A(\cdot) = \{A(t); t \ge 0\}$ which has bounded variations. For example, the $Z(\cdot)$ of \eq{Z-SIE} is a semi-martingale if it exists. We assume that the martingale $M(\cdot)$ has finite quadratic variations $[X]_{\bullet} \equiv \{[X]_{t}; t \ge 0\}$. For this semi-martingale $X(\cdot)$, local time $L_{x}(t)$ for $x \in \dd{R}$ and $t \ge 0$ is defined through
\begin{align}
\label{eq:L-time1}
  \int_{-\infty}^{\infty} L_{x}(t) p(x) dx = \int_{0}^{t} p(X(s)) d[X]_{s} \quad \mbox{ for any measurable function $p$}.
\end{align}
See Theorem 7.1 of \cite{KaraShre1998} for details about the definition of local time. Note that the local time of \cite{KaraShre1998} is half of the local time in this paper. Applying $p(y) = 1_{(x-\varepsilon,x+\varepsilon)}(y)$ for $\varepsilon > 0$ to \eq{L-time1}, we can see that
\begin{align}
\label{eq:L-time2}
  L_{x}(t) = \lim_{\varepsilon \downarrow 0} \frac 1{2\varepsilon} \int_{0}^{t} 1_{(x-\varepsilon,x+\varepsilon)}(X(s)) d[X]_{s}, \quad a.s. \quad x \in \dd{R}, t \ge 0
\end{align}
This can be used as the definition of the local time. 

There are two versions of the local time since $L_{x}(t)$ is continuous in $t$, but may not be continuous in $x$. So, usually, the local time $L_{x}(t)$ is assumed to be right-continuous. However, if the finite variation component of $X(\cdot)$ is not atomic, then $L_{x}(t)$ is continuous in $x$ (see Theorem 22.4 of \cite{Kall2001}). Hence,

\begin{lemma}
\label{lem:LT-time}
For the $Z(\cdot)$ of a state-dependent reflected diffusion on $[0,\infty)$, its local time $L_{x}(t)$ is continuous in $x$ for each $t \ge 0$. Furthermore, $\dd{E}[L_{x}(t)]$ is finite by \eq{L-time1} for $X(\cdot) = Z(\cdot)$.
\end{lemma}
 
\subsection{Generalized Ito formula}
\label{app:G-Ito}

We next introduce a generalized Ito formula for a semi-martingale $X(\cdot)$ and a convex test function $f$, which is given by
\begin{align}
\label{eq:G-Ito-X}
  f(X(t)) & = f(X(0)) + \int_{0}^{t} f'(X(s)-) dX(s) + \frac 12 \int_{0}^{\infty} L_{u}(t) \mu_{f}(du), \qquad t \ge 0,
\end{align}
where $L_{u}(t)$ is the local time of $X(\cdot)$ which is right-continuous in $u \in \dd{R}$, and $\mu_{f}$ on $(\dd{R}_{+}, \sr{B}(\dd{R}_{+}))$ is a measure on $(\dd{R},\sr{B}(\dd{R})$, defined by
\begin{align}
\label{eq:Mf-1}
  \mu_{f}([u,v)) = f'(v-) - f'(u-), \qquad u < v \mbox{ in } \dd{R},
\end{align}
where  $f'(u-)$ is the left derivative of $f$ at $u$. The formula \eq{G-Ito-X} is also called Ito-Meyer-Tanaka formula (e.g., see Theorem 6.22 of \cite{KaraShre1998} and Theorem 22.5 of \cite{Kall2001}). 

Note that if $f$ is a concave function, then $-f$ is convex, so \eq{G-Ito-X} yields
\begin{align}
\label{eq:G-Ito-X-concave}
  f(X(t)) & = f(X(0)) + \int_{0}^{t} f'(X(s)-) dX(s) - \frac 12 \int_{0}^{\infty} L_{u}(t) \mu_{(-f)}(du), \qquad t \ge 0,
\end{align}

The excellence of \eq{G-Ito-X} is that test function $f$ does not need to be continuously twice differentiable. In particular, if $f'(y)$ is almost surely differentiable for $y \in [u,v)$, then
\begin{align}
\label{eq:Mf-2}
  \mu_{f}([u,v)) = \int_{u}^{v} f''(y-) dy, \qquad u < v \mbox{ in } \dd{R},
\end{align}
Hence, if $f$ and $f'$ are almost surely differentiable with respect to Lebesgue measure, then, applying \eq{L-time1}, \eq{G-Ito-X} can be written as
\begin{align}
\label{eq:Ito1}
  f(X(t)) = f(X(0)) + \int_{0}^{t} f'(X(s)) dX(s) + \frac 12 \int_{0}^{\infty} f''(X(s)) d[X]_{s}, \qquad t \ge 0.
\end{align}
This formula is also called a generalized version of Ito formula because test function $f$ does not need to be continuously twice differentiable.

We apply the generalized Ito formula \eq{G-Ito-X} to the solution of $Z(t)$ of \eq{Z-SIE}. For this, define $M(\cdot) \equiv \{M(t); t \ge 0\}$ by
\begin{align}
\label{eq:Z-M}
  M(t) \equiv \int_{0}^{t} \sigma(X(s)) dW(s), \qquad t \ge 0,
\end{align}
then $M(\cdot)$ is a martingale. Denote the quadratic variations of $Z(\cdot)$ and $M(\cdot)$, respectively, by $\qvs{Z}_{t}$ and $\qvs{M}_{t}$. Since $Z(t)$ and $Y(t)$ are continuous in $t$ (see \rem{Z-X-existence}), it follows from \eq{Z-SIE} that
\begin{align}
\label{eq:M-qv}
  \qvs{Z}_{t} & = \qvs{M}_{t} = \int_{0}^{t} \sigma^{2}(Z(s)) ds,  \qquad t \ge 0.
\end{align}

Hence, the generalized Ito formula \eq{G-Ito-X} can be written as
\begin{align}
\label{eq:G-Ito-Z}
  f(Z(t)) & = f(Z(0)) + \int_{0}^{t} f'(Z(s)-) b(Z(s)) ds +  \int_{0}^{t} f'(Z(s)-) \sigma(Z(s)) dW(s)  \nonumber\\
  & \quad + \frac 12 \int_{0}^{\infty} L_{u}(t) \mu_{f}(du) + \int_{0}^{t} f'(Z(s)) Y(ds), \qquad t \ge 0.
\end{align}

\subsection{Tanaka formula}
\label{app:Tanaka-formula}

For constant $a \in \dd{R}$, let $f(x) = |x-a|$ for \eq{G-Ito-X}, then $f'(x) = 1(x > a) - 1(x < a)$ for $x \not= a$. Hence, $\mu_{f}(B) = 21(a \in B)$, and thereforet follows from \eq{G-Ito-X} that
\begin{align}
\label{eq:Tanaka}
 & |X(t)-a| = |X(0)-a| + \int_{0}^{t} \sgn(X(s)- a) dX(s) + L_{a}(t), \qquad t \ge 0,
\end{align}
where $\sgn(x) = 1(x>0) - 1(x \le 0)$. This formula is called a Tanaka formula because it is originally studied for Brownian motion by \citet{Tana1963}. 

\subsection*{Acknowledgements}

This paper partly follows up my talk at Uma Prabhu centennial conference in September, 2024.  I am grateful to its organizers. I also benefited from discussions with Jim Dai.

%\bibliography{../../../texmf/bib/dai20230708M2}
%\end{document}

\end{document}